\documentclass{article}
\usepackage{amsmath}
\usepackage{amssymb}
\usepackage{amsfonts}
\usepackage{tikz}
\usepackage{color}
\usepackage{graphicx,caption}

\def\gcd#1#2{\text{gcd}(#1,#2)}
\def\S#1#2{{\cal F}_{#1}^{#2}}
\def\N{\mathbb N}
\def\pi{\text{id}}

\newcounter{NN}
\begin{document}
\bibliographystyle{plain}
\title{On the Fourier transform of \\ the greatest common divisor} 
\author{Peter H.~van der Kamp}
\date{Department of Mathematics and Statistics\\
La Trobe University \\
Victoria 3086, Australia\\[5mm]
\today
}

\pagestyle{plain} \maketitle

\begin{abstract}
The discrete Fourier transform of the greatest common divisor
\[
\widehat{\pi}[a](m)=\sum_{k=1}^m \gcd{k}{m} \alpha_m^{ka},
\]
with $\alpha_m$ a primitive $m$-th root of unity, is a multiplicative
function that generalises both the gcd-sum function and Euler's totient function.
On the one hand it is the Dirichlet convolution of the identity with Ramanujan's sum, %\cite{WS},
$\widehat{\pi}[a]=\pi\ast c_\bullet(a)$, and on the other hand it can be written as a generalised
convolution product, %\cite{AA}
$\widehat{\pi}[a]=\pi \ast_a \phi$. 

We show that $\widehat{\pi}[a](m)$ counts the number of elements in the set of ordered pairs $(i,j)$ such
that $i\cdot j \equiv a \mod m$. Furthermore we generalise a dozen known identities for the totient function,
to identities which involve the discrete Fourier transform of the greatest common divisor, including its partial
sums, and its Lambert series.
\end{abstract}

\section{Introduction}
In \cite{WS} discrete Fourier transforms of functions of the greatest common
divisor were studied, i.e. 
\[
\widehat{h}[a](m) = \sum_{k=1}^m h(\gcd{k}{m})\alpha_m^{ka},
\]
where $\alpha_m$ is a primitive $m$-th root of unity. The main result in that
paper is the identity\footnote{Similar results in the context of $r$-even function were
obtained earlier, see \cite{LTPH} for details.} $\widehat{h}[a]=h\ast c_\bullet(a)$, where $\ast$ denoted Dirichlet
convolution, i.e.
\begin{equation} \label{hffc}
\widehat{h}[a](m) = \sum_{d\mid m} h(\frac{m}{d}) c_d(a),
\end{equation}
and
\begin{equation} \label{Ram}
c_m(a)=\sum_{\underset{\gcd{k}{m}=1}{k=1}}^m \alpha_m^{ka} 
\end{equation}
denotes Ramanujan's sum. Ramanujan's sum generalises both Euler's totient function
$\phi=c_\bullet(0)$ and the M\"obius function $\mu=c_\bullet(1)$. Thus, identity (\ref{hffc})
generalizes the formula
\begin{equation} \label{C}
\sum_{k=1}^m h(\gcd{k}{m}) = (h \ast \phi)(m),
\end{equation}
already known to Ces\`aro in 1885. The formula (\ref{hffc}) shows that
$\widehat{h}[a]$ is multiplicative if $h$ is multiplicative (because $c_\bullet(a)$ is multiplicative
and Dirichlet convolution preserves multiplicativity).

\def\g{{\mathcal P}}
Here we will take $h=\pi$ to be the identity function (of the greatest common divisor) and study
its Fourier transform.
Obviously, as $\pi: \pi(n)=n$ is multiplicative, the function $\widehat{\pi}[a]$ is multiplicative, for all $a$. 
Two special cases are well-known. Taking $a=0$ we have $\widehat{\pi}[0]=\g$, where
\begin{equation} \label{g}
\g(m)=\sum_{k=1}^m \gcd{k}{m}.
\end{equation}
is known as Pillai's arithmetical function or the gcd-sum function.
Secondly, by taking $a=1$ in (\ref{hffc}), we find that $\widehat{\pi}[1]=\pi \ast \mu$ equals $\phi$,
by M\"obius inversion of Euler's identity $\phi \ast u = \pi$, where $u=\mu^{-1}$ is the unit function defined by $u(m)=1$.

Let $\S am$ denote the set of ordered pairs $(i,j)$ such that $i\cdot j\equiv a \mod m$,
the set of {\em factorizations of $a$ modulo $m$}.
We claim that $\widehat{\pi}[a](m)$ counts its number of elements. Let us consider the mentioned special cases.
\begin{itemize}
\item[$a=0$] For given $i\in\{1,2,\ldots,m\}$ the congruence $i\cdot j \equiv 0 \mod m$ yields
\[
\frac{i}{\gcd im} j \equiv 0 \mod  \frac{m}{\gcd im},
\]
which has a unique solution modulo $m/\gcd im$, and so there are $\gcd im$ solutions
modulo $m$. Hence, the total number of elements in $\S 0m$ is $\g(m)$.
\item[$a=1$] The totient function $\phi(m)$ counts the number of invertible congruence classes
modulo $m$. As for every invertible congruence class $i$ modulo $m$ there is a unique
$j=i^{-1} \mod m$ such that $i\cdot j \equiv 1 \mod m$, it counts the number of
elements in the set $\S 1m$.
\end{itemize}

\noindent
To prove the general case we employ a Kluyver-like formula for $\widehat{\pi}[a]$, that is, a formula similar
to the formula for the Ramanujan sum function
\begin{equation} \label{klu}
c_k(a)=\sum_{d\mid \gcd ak} d \mu(\frac{k}{d}).
\end{equation}
attributed to Kluyver. Together the identities (\ref{hffc}) and (\ref{klu}) imply, cf. section \ref{grs},
\begin{equation} \label{pham}
\widehat{\pi}[a](m)=\sum_{d\mid \gcd am} d \phi(\frac{m}{d}),
\end{equation}
and we will show, in the next section, that the number of factorizations of $a \mod m$ is
given by the same sum.

The right hand sides of (\ref{klu}) and (\ref{pham}) are particular instances of the so called generalized
Ramanujan sums \cite{AA}, and both formulas follow as consequence of a general formula for the Fourier coefficients
of these generalised Ramanujan sums \cite{Tom,Tom2}.
In section \ref{grs} we provide simple proofs for some of the nice properties of these sums. In particular
we interpret the sums as a generalization of Dirichlet convolution. This interpretation lies at the heart of many
of the generalised totient identities we establish in section \ref{gti}.

\section{The number of factorizations of $a \mod m$} \label{numset}
For given $i,m\in\N$, denote $g=\gcd im$. If the congruence $i\cdot j \equiv a \mod m$ has a solution
$j$, then $g\mid a$ and $j\equiv i^{-1} a/g$ is unique mod $m/g$, so mod $m$ there are $g$ solutions.
This yields
\[
\#\S am=\sum_{\underset{ \gcd im\mid a}{i=1}}^m \gcd im,
\]
which can be written as
\begin{equation} \label{f}
\#\S am=\sum_{d\mid a} \sum_{\underset{\gcd im=d}{i=1}}^m d
\end{equation}
If $d\nmid m$ then the sum \[
\sum_{\underset{\gcd im=d}{i=1}}^m 1\]
is empty. Now let $d\mid m$. The only integers $i$ which contribute to the sum
are the multiples of $d$, $kd$, where $\gcd k{m/d}=1$. There are exactly $\phi(m/d)$ of them.
Therefore the right hand sides of formulae (\ref{pham}) and (\ref{f}) agree, and hence $\#\S am=\widehat{\pi}[a](m)$.

\section{A historical remark, and generalised Ramanujan sums} \label{grs}
It is well known that Ramanujan was not the first who considered the sum $c_m(a)$.
Kluyver proved his formula (\ref{klu}) in 1906, twelve years before Ramanujan
published the novel idea of expressing arithmetical functions in the form of a series
$\sum_s a_s c_s(n)$ \cite{Ram}. It is not well known that Kluyver actually showed that
$c_m(a)$ equals Von Sterneck's function, introduced in \cite{VoS}, i.e.
\begin{equation} \label{Hol}
c_m(a)=\frac{\mu(\frac{m}{\gcd am})\phi(m)}{\phi(\frac{m}{\gcd am})}.
\end{equation}
This relation is referred to in the literature as H\"older's relation, cf. the remark on page 213 in \cite{AA}. However, H\"older published
this relation thirty years after Kluyver \cite{Hol}. We refer to \cite[Theorem 2]{AA}, or \cite[Theorem 8.8]{Tom}
for a generalisation of (\ref{Hol}). The so called generalized Ramanujan sums,
\begin{equation} \label{fag}
f\ast_a g(m)=\sum_{d\mid\gcd am} f(d)g(\frac{m}{d}),
\end{equation}
were introduced in \cite{AA}. The notation $\ast_a$ is new, the sums are denoted $S(a;m)$ in \cite{AA}, $s_m(a)$ in
\cite{Tom}, and $S_{f,h}(a,m)$ in \cite{Tom2}. In the context of $r$-even functions \cite{LTPH} the sums are denoted $S_{f,g}(a)$, and
considered as sequences of $m$-even functions, with argument $a$. We consider the sums as a sequence of functions with
argument $m$, labeled by $a$. We call $f\ast_a g$ {\bf the ${\mathbf a}$-convolution of ${\mathbf f}$ and ${\mathbf g}$}.

The concept of $a$-convolution is a generalization of Dirichlet convolution as $f\ast_0 g=f\ast g$. As we will see below, the function $f\ast_a g$ is
multiplicative (for all $a$) if $f$ and $g$ are, and the following inter-associative property holds, cf. \cite[Theorem 4]{Tom2}. 
\begin{equation} \label{assp}
(f\ast_a g) \ast h = f\ast_a (g \ast h).
\end{equation}
We also adopt the notation $f_a=\pi \ast_a f$, and call this {\bf the  Kluyver, or ${\mathbf a}$-extension of ${\mathbf f}$}. Thus, we have $f_0=\pi \ast f$, $f_1=f$, and formulas
(\ref{klu}) and (\ref{pham}) become $c_m(a)=\mu_a(m)$, and $\widehat{\pi}[a]=\phi_a$, respectively.

The identity function $I$, defined by $f\ast I=f$, is given by $I(k)=[k=1]$, where the Iverson bracket is, with $P$ a logical statement,
\[
[P]=\left\{\begin{array}{ll} 1 & \text{if } P, \\ 0 & \text{if } \text{not } P. \end{array} \right.
\]
Let us consider the function $f\ast_a I$. It is
\[
f \ast_a I(k)=\sum_{d\mid\gcd ak} f(d)[d=k]=[k\mid a]f(k).
\]
Since the function $k\rightarrow[k\mid a]$ is multiplicative, the function $f\ast_a I$ is multiplicative if $f$ is multiplicative. Also, we may write, cf. \cite[eq. (9)]{Tom2},
\[
f\ast_a g(m)=\sum_{d\mid m} [d\mid a]f(d)g(\frac{m}{d}) = (f\ast_a I) \ast g (m),
\]
which shows that $f\ast_a g$ is multiplicative if $f$ and $g$ are. Also, the inter-associativity property (\ref{assp}) now easily follows from
the associativity of the Dirichlet convolution,
\[
 (f\ast_a g) \ast h = ((f\ast_a I) \ast g) \ast h = (f\ast_a I) \ast (g \ast h) = f\ast_a (g \ast h).
\]
We note that the $a$-convolution product is neither associative, nor commutative. The inter-associativity and the commutativity of
Dirichlet convolution imply that
\begin{align} \label{ig}
f_a \ast g = (f\ast g)_a  = f \ast g_a.
\end{align}

Formula (\ref{pham}) states that the Fourier transform of the greatest common divisor is the Kluyver extension
of the totient function. We provide a simple proof.

\smallskip
\noindent
{\bf Proof} [of (\ref{pham})]  Employing (\ref{hffc}), (\ref{klu}) and (\ref{ig}) we have
$
\widehat{\pi}[a] = \pi \ast c_\bullet(a) = \pi \ast \mu_a = (\pi \ast \mu)_a = \phi_a.
$
\hfill $\square$

\smallskip
\noindent
The formula (\ref{pham}) also follows as a special case of the following formula for the Fourier coefficients of
$a$-convolutions, 
\begin{equation} \label{ha}
f\ast_a g(m)=\sum_{k=1}^m h_k(m) \alpha_m^{ka},\qquad  h_k = g \ast_k \frac{f}{\pi},
\end{equation}
given in \cite{AA,Tom}. The formula (\ref{ha}) combines with (\ref{hffc}) and (\ref{klu}) to yield a formula for functions
of the greatest common divisor, $\bar{h}[k]: m\rightarrow h(\gcd km)$, namely
\begin{equation} \label{con}
\bar{h}[k] = (h \ast \mu) \ast_k u.
\end{equation}
\noindent
{\bf Proof} [of (\ref{con})] The Fourier coefficients of $\widehat{h}[a](m)$ are $\bar{h}[k](m)$. But $\widehat{h}[a]=h \ast c_\bullet(a)
= (\pi \ast_a \mu) \ast h = \pi \ast_a (\mu \ast h)$, and so, using  (\ref{ha}), the Fourier coefficients are also given by $(h \ast \mu) \ast_k u(m)$.
\hfill $\square$

\smallskip
\noindent
For a Dirichlet convolution with a Fourier transform of a function of the greatest common divisor we have
\begin{equation} \label{mt1}
f\ast\widehat{g}[a] = \widehat{f\ast g}[a].
\end{equation}
\noindent
{\bf Proof} [of (\ref{mt1})] $f\ast\widehat{g}[a]=f\ast(g\ast\mu_a)=(f\ast g)\ast \mu_a= \widehat{f\ast g}[a]$ \hfill $\square$

\smallskip
\noindent
Similarly, for an $a$-convolution with a Fourier transform of a function of the greatest common divisor, 
\begin{equation} \label{mt2}
f\ast_a \widehat{g}[b] = \widehat{f\ast_a g}[b].
\end{equation}
\noindent
{\bf Proof} [of (\ref{mt2})] $f\ast_a\widehat{g}[b]=f\ast_a(g\ast\mu_b)=(f\ast_a g)\ast \mu_b= \widehat{f\ast_a g}[b]$ \hfill $\square$

\section{Generalised totient identities} \label{gti}
The totient function is an important function in number theory, and related fields of mathematics. It is extensively studied, connected to many other notions and functions, and there exist numerous generalisation and extensions, cf. the chapter "The many facets of Euler's totient" in \cite{SC}. The Kluyver extension of the totient function is a very natural extension, and it is most surprising it has not been studied before.
In this section we generalise a dozen known identities for the totient function $\phi$, to identities which involve its Kluyver extension $\phi_a$,
a.k.a. the discrete Fourier transform of the greatest common divisor. This includes a generalisation of Euler's identity, the partial sums of $\phi_a$, and its Lambert series.

\subsection{The value of $\phi_a$ at powers of primes}
We start by providing a formula for the value of $\phi_a$ at powers of primes. This depends only on the multiplicity
of the prime in $a$. The formulae, with $p$ prime,
\[
\g(p^k) = (k + 1)p^k - kp^{k-1}, \qquad \phi(p^k)=p^k-p^{k-1},
\]
of which the first one is Theorem 2.2 in \cite{Bro}, generalise to
\begin{equation} \label{pp}
\phi_a(p^k)=\left\{\begin{array}{ll}
(p^k-p^{k-1})(l+1) & \text{if } l<k, \\
(k+1)p^{k}-kp^{k-1} & \text{if } l\geq k,
\end{array}
\right.
\end{equation}
where $l$ is the largest integer, or infinity, such that $p^l\mid a$.

\smallskip
\noindent
{\bf Proof} [of (\ref{pp})] We have
\begin{align*}
\phi_a(p^k)&=\sum_{d\mid \gcd {p^l}{p^k}} d \phi(\frac{p^k}{d})\notag\\
&=\sum_{r=0}^{\min(l,k)} p^r \phi(p^{k-r})\notag \\
&=\left\{\begin{array}{ll}
\sum_{r=0}^l p^k-p^{k-1} & \text{if } l<k, \\
(\sum_{r=0}^{k-1} p^k-p^{k-1}) +  p^k & \text{if } l\geq k,
\end{array}
\right.
\end{align*}
which equals (\ref{pp}).

\subsection{Partial sums of $\phi_a/\pi$}
To generalise the totient identity 
\begin{equation} \label{fti}
\sum_{k=1}^n \frac{\phi(k)}{k} = \sum_{k=1}^n \frac{\mu(k)}{k} \lfloor \frac{n}{k} \rfloor.
\end{equation}
to an identity for $\phi_a$ we first establish
\begin{equation} \label{ook}
\sum_{k=1}^n \frac{f_0(k)}{k} = \sum_{k=1}^n \frac{f(k)}{k} \lfloor \frac{n}{k} \rfloor.
\end{equation}

\smallskip
\noindent
{\bf Proof} [of (\ref{ook})] Since there are $\lfloor n/d \rfloor$ multiples of $d$ in the range $[1,n]$ it follows that
\[
\sum_{k=1}^n \frac{f\ast \pi(k)}{k} = \sum_{k=1}^n \sum_{d\mid k} \frac{f(d)}{d} = \sum_{d=1}^n \frac{f(d)}{d} \lfloor \frac{n}{d} \rfloor.
\]\hfill $\square$

As a corollary we obtain
\begin{equation} \label{cook}
\sum_{k=1}^n \frac{f\ast_a g_0(k)}{k} = \sum_{k=1}^n \frac{f\ast_a g(k)}{k} \lfloor \frac{n}{k} \rfloor.
\end{equation}

\smallskip
\noindent
{\bf Proof} [of (\ref{cook})] Employing (\ref{assp}) we find
\[
\sum_{k=1}^n \frac{f\ast_a (g \ast \pi)(k)}{k} = \sum_{k=1}^n \frac{(f\ast_a g) \ast \pi(k)}{k} = \sum_{k=1}^n \frac{f\ast_a g(k)}{k} \lfloor \frac{n}{k} \rfloor.
\]\hfill $\square$

Now taking $f=\pi$ and $g=\mu$ in (\ref{cook}) we find
\begin{equation} \label{g1}
\sum_{k=1}^n \frac{\phi_a(k)}{k} = \sum_{k=1} ^n \frac{c_k(a)}{k} \lfloor \frac{n}{k} \rfloor.
\end{equation}

\subsection{Partial sums of $\g_a/\pi$ expressed in terms of $\phi_a$}
Taking $f=\pi$ and $g=\phi$ in (\ref{cook}) we find
\begin{equation} \label{g2}
\sum_{k=1}^n \frac{\g_a(k)}{k} = \sum_{k=1} ^n \frac{\phi_a(k)}{k} \lfloor \frac{n}{k} \rfloor.
\end{equation}

Note that by taking either $a=0$ in (\ref{g1}), or $a=1$ in (\ref{g2}), we find an identity involving
the totient function and the gcd-sum function,
\begin{equation} \label{g3}
\sum_{k=1}^n \frac{\g(k)}{k} = \sum_{k=1} ^n \frac{\phi(k)}{k} \lfloor \frac{n}{k} \rfloor.
\end{equation}

\subsection{Partial sums of $\phi_a$}
To generalise the totient identity, with $n>0$,
\begin{equation} \label{iop}
\sum_{k=1}^n \phi(k) = \frac{1}{2}\left(1+\sum_{k=1}^n \mu(k) \lfloor \frac{n}{k} \rfloor^2\right),
\end{equation}
we first establish
\begin{equation} \label{fid}
\sum_{k=1}^n f_0 (k) = \frac{1}{2}\left( \sum_{k=1}^n f(k) \lfloor \frac{n}{k} \rfloor^2 + \sum_{k=1}^n f\ast u (k)\right).
\end{equation}

\smallskip
\noindent
{\bf Proof} [of (\ref{fid})] We have, by changing variable $k=dl$,
\begin{align*}
\sum_{k=1}^n (2 f\ast \pi - f\ast u)(k) &= \sum_{k=1}^n \sum_{d\mid k} f(d) (\frac{2k}{d} -1) \\
&= \sum_{d=1}^n \sum_{l=1}^{\lfloor n/d \rfloor} f(d) (2l -1) \\
&= \sum_{d=1}^n f(d) \lfloor \frac{n}{d} \rfloor^2.
\end{align*}\hfill $\square$

\newpage
\noindent
Note that this gives a nice proof of (\ref{iop}), taking $f=\mu$, as $\sum_{k=1}^n I(k)=[k>0]$. When $f=\mu_a$, then (\ref{ig}) implies
$f\ast \pi=\phi_a$, and $f\ast u= I_a$, and therefore as a special case of (\ref{fid}) we obtain
\begin{equation} \label{g4}
\sum_{k=1}^n \phi_a(k) = \frac{1}{2}\left( \sum_{k\mid a} k[k\leq n] +\sum_{k=1}^n c_k(a) \lfloor \frac{n}{k} \rfloor^2\right).
\end{equation}
We remark that when $n\geq a$ we have
$
\sum_{k\mid a} k[k\leq n]  = \sigma(a),
$
where $\sigma=\pi\ast u$ is the sum of divisors function.

%We note that the expression (\ref{rhs}) counts the number $n_d$ of lattice points $(i,j)$ in the square $[1,n]\times[1,n]$,
%such that $gcd(i,j)=d$ and computes the sum $\sum_{d\mid a} d n_d$. The change that if one pick 

\subsection{Generalisation of Euler's identity}%\smallskip
Euler's identity, $\phi\ast u=\pi$, generalises to
\begin{equation} \label{jkl}
\sum_{d\mid m} \phi_a(d) = \tau(\gcd am)m,
\end{equation}
where $\tau=u\ast u$ is the number of divisors function.

\smallskip
\noindent
{\bf Proof} [of (\ref{jkl})] We have $\phi_a\ast u=(\phi\ast u)_a=\pi_a$ where
\begin{equation} \label{ida}
\pi_a(m) = \sum_{d\mid \gcd am} d \frac{m}{d} = m \tau(\gcd am).
\end{equation} \hfill $\square$

%\smallskip
\subsection{Partial sums of $\g_a$ expressed in terms of $\phi_a$ (and $\tau$)}
If $f=\phi_a$, then $f\ast \pi=\g_a$, and (\ref{fid}) becomes, using (\ref{jkl}),
\begin{equation} \label{g5}
\sum_{k=1}^n \g_a(k) = \frac{1}{2}\left( \sum_{k=1}^n \tau(\gcd ak)k +\sum_{k=1}^n \phi_a(k) \lfloor \frac{n}{k} \rfloor^2\right).
\end{equation}

%\smallskip
\subsection{Four identities of C\'esaro}
According to Dickson \cite{Dick} the following three identities were obtained by C\'esaro:
\begin{align}
\sum_{d\mid n} d \phi(\frac{n}{d})&=\g(n), \label{c1}\\
\sum_{d\mid n} \frac{d}{n} \phi(d)&=\sum_{j=1}^n \frac{1}{\gcd jn}, \label{c2}\\
\sum_{d\mid n} \phi(d)\phi(\frac{n}{d})&=\sum_{j=1}^n\phi(\gcd jn). \label{c3}\
\end{align}
Identity (\ref{c1}), which is Theorem 2.3 in \cite{Bro}, is obtained by taking $a=0$ in (\ref{pham}), or $h=\pi$ in (\ref{C}).
It generalises to
\begin{equation} \label{p1}
\sum_{d\mid n} d\phi_a(\frac{n}{d}) = \g_a(n).
\end{equation}

\smallskip
\noindent
{\bf Proof} [of (\ref{p1})] By taking $f=\phi$ and $g=\pi$ in (\ref{ig}). \hfill $\square$

\smallskip
\noindent 
Identity (\ref{c2}) is obtained by taking $h=1/\pi$ in (\ref{C}) and generalises to
\begin{equation} \label{p2}
\sum_{d\mid n} \frac{d}{n}\phi_a(d) = \sum_{j=1}^n \sum_{d\mid\gcd an} \frac{1}{\gcd j{\frac{n}{d}}},
\end{equation}

\smallskip
\noindent
{\bf Proof} [of (\ref{p2})] By taking $f=\phi$ and $g=1/\pi$ in (\ref{ig}). \hfill $\square$

\smallskip
\noindent
Identity (\ref{c3}) is also a special case of (\ref{C}), with $h=\phi$. It generalises to
\begin{equation} \label{p3}
\sum_{d\mid n} \phi_a(d)\phi_b(\frac{n}{d})=\sum_{j=1}^n \sum_{d\mid\gcd an} \phi_b(\gcd j{\frac{n}{d}}).
\end{equation}

\smallskip
\noindent
{\bf Proof} [of (\ref{p3})] We have
\[
(\pi \ast_a \phi) \ast (\pi \ast_b \phi)  = \pi \ast_a (\pi \ast_b (\phi\ast \phi))
= \pi \ast_a (\pi \ast_b \widehat{\phi}[0]) 
= \pi \ast_a \widehat{\phi_b}[0],
\]
and evaluation at $m$ yields
\[
\sum_{d\mid\gcd am} d \sum_{j=1}^{m/d} \phi_b(\gcd j{\frac{m}{d}}) =
\sum_{d\mid\gcd am} \sum_{j=1}^{m} \phi_b(\gcd j{\frac{m}{d}}).
\]
\hfill $\square$

\smallskip
\noindent
The more general identity (\ref{C}) generalises to
\begin{equation} \label{gC}
\sum_{k=1}^m h_a(\gcd km) = h\ast \phi_a (m).
\end{equation}

\subsection{Three identities of Liouville}
Dickson \cite[p.285-286]{Dick} states, amongst many others identities that were presented by Liouville in the series \cite{Liou}, the following
\begin{align}
\sum_{d\mid m} \phi(d) \tau(\frac{m}{d}) = \sigma(m), \label{l1}\\
\sum_{d\mid m} \phi(d) \sigma[n+1](\frac{m}{d}) = m\sigma[n](m), \label{l2}\\ 
\sum_{d\mid m} \phi(d) \tau(\frac{m^2}{d^2}) = \sum_{d\mid m} d\theta(\frac{m}{d}), \label{l3}
\end{align}
where $\sigma[n]=\pi[n]\ast u$, $\pi[n](m)=m^n$, and $\theta(m)$ is the number of decompositions of $m$ into two relatively prime
factors. All three are of the form $\phi\ast f = g$ and therefore they gain significance due to (\ref{C}), thought
Liouville might not have been aware of this. For example, (\ref{C}) and (\ref{l1}) combine to yield
\[
\sum_{k=1}^m \tau(\gcd km) = \sigma(m).
\]
The three identities are easily proven by substituting $\tau=u\ast u$, $\sigma[n]=\pi[n]\ast u$, $\tau\circ\pi[2]=\theta\ast u$, $\phi=\mu\ast \pi$,
and using $\mu\ast u=I$. They generalise to
\begin{align}
\sum_{d\mid m} \phi_a(d) \tau(\frac{m}{d}) = \sigma_a(m), \label{lg1}\\
\sum_{d\mid m} \phi_a(d) \sigma[n+1](\frac{m}{d}) = m u \ast_a \sigma[n](m), \label{lg2}\\ 
\sum_{d\mid m} \phi_a(d) \tau(\frac{m^2}{d^2}) = \sum_{d\mid m} d\tau(\gcd ad)\theta(\frac{m}{d}). \label{lg3}
\end{align}
These generalisation are proven using the same substitutions, together with (\ref{ig}), or for the latter identity, (\ref{assp}) and (\ref{ida}).
\subsection{One identity of Dirichlet}\smallskip
Dickson \cite{Dick} writes that Dirichlet \cite{Dir}, by taking partial sums on both sides of Euler's identity, obtained
\[
\sum_{k=1}^n \lfloor \frac{n}{k} \rfloor \phi(k) = {n+1\choose 2}.
\]
By taking partial sums on both sides of equation (\ref{jkl}) we obtain
\begin{equation} \label{tyu}
\sum_{k=1}^n \lfloor \frac{n}{k} \rfloor \phi_a(k) = \sum_{d\mid a} d {\lfloor \frac{n}{d} \rfloor +1\choose 2} .
\end{equation}

\smallskip
\noindent
{\bf Proof} [of (\ref{tyu})]  Summing the left hand side of (\ref{jkl}) over $m$ yields
\[
\sum_{m=1}^n \sum_{d \mid m} \phi_a(d) = \sum_{d=1}^n \lfloor \frac{n}{d} \rfloor  \phi_a(d)
\]
and summing the right hand side of (\ref{jkl}) over $m$ yields
%\begin{align*}
\[
\sum_{m=1}^n \tau(\gcd am)m = \sum_{m=1}^n \sum_{d\mid\gcd am} m  
= \sum_{d\mid a}  \sum_{k=1}^{\lfloor n/d\rfloor} dk 
= \sum_{d\mid a} d\lfloor \frac{n}{d}\rfloor\left(\lfloor  \frac{n}{d}\rfloor+1\right)/2.
\]%\end{align*}
 \hfill $\square$

\subsection{The Lambert series of $\phi_a$}
As shown by Liouville \cite{Liou}, cf. \cite[p.120]{Dick}, the Lambert series of the totient function is given by
\[
\sum_{m=1}^\infty \phi(m) \frac{x^m}{1-x^m} = \frac{x}{(1-x)^2}.
\]
The Lambert series for $\phi_a$ is given by
\begin{equation} \label{lam}
\sum_{m=1}^\infty \phi_a(m) \frac{x^m}{1-x^m} = p[a](x)\frac{x}{(1-x^a)^2},
\end{equation}
where the coefficients of $p[a](x)=\sum_{k=1}^{2a} c[a](k)x^{k-1}$ are given by
$c[a]=\pi_a \circ t[a]$, and $t[a]$ is the piece-wise linear function $t[a](n)=a-|n-a|$.

The polynomials $p[a]$ seem to be irreducible over ${\mathbb Z}$ and their zeros are in some sense
close to the $a$-th roots of unity, see Figures \ref{37} and \ref{35}.

\noindent
\begin{center}
\includegraphics[width=6cm]{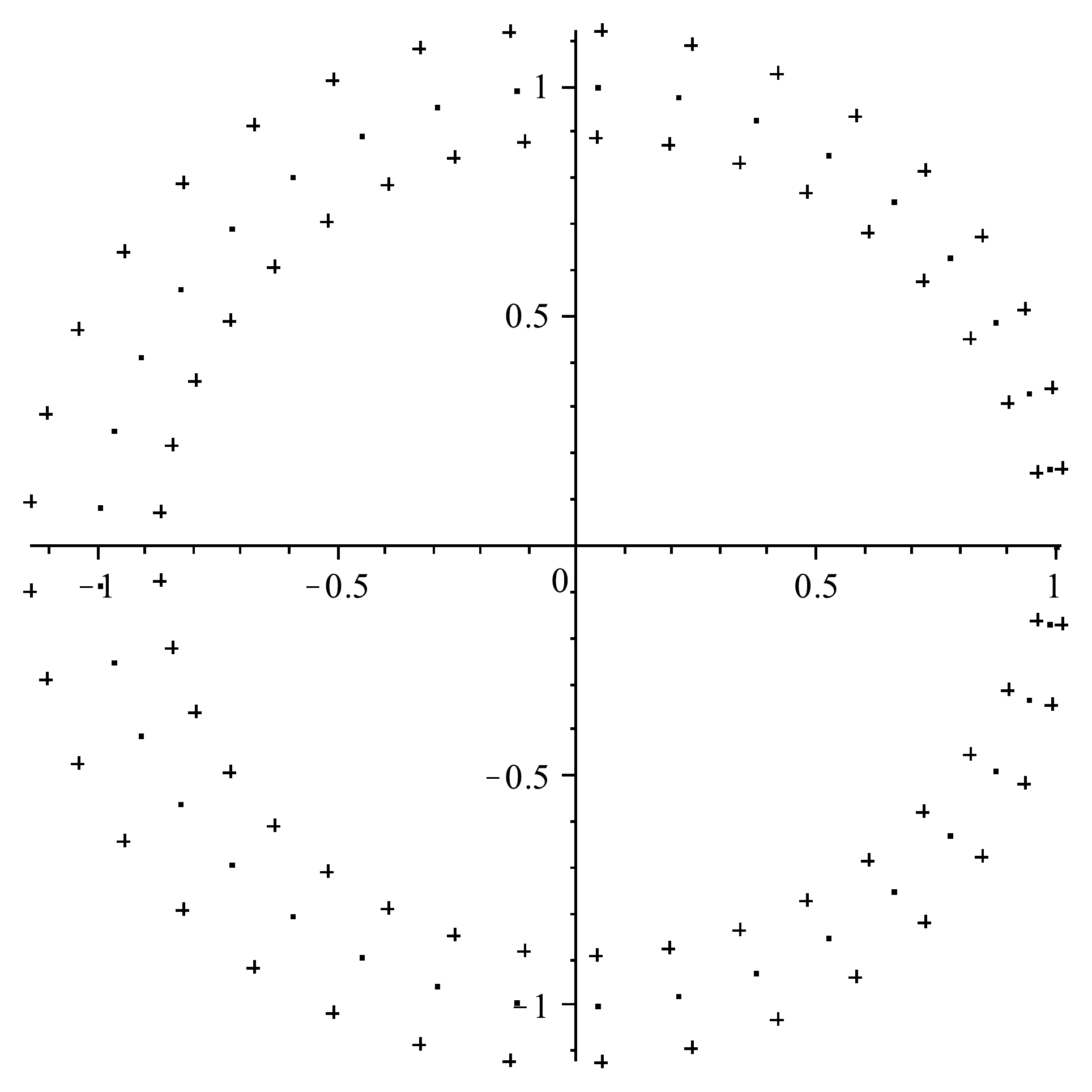}
\captionof{figure}{\label{37} The roots of $p[37]$ are depicted as crosses and the $37^{\text{th}}$ roots of unity as points.
This figure shows that when $a$ is prime the roots of $p[a]$ that are close to 1 are closer to $a^{\text{th}}$ roots of unity.}
\end{center}

\smallskip
\noindent
{\bf Proof} [of (\ref{lam})] Ces\`aro proved \cite{Dick}
\[
 \sum_{n=1}^\infty f(n)\frac{x^n}{1-x^n} = \sum_{n=1}^\infty x^n \sum_{d\mid n} f(d),
\]
cf. exercise 31 to chapter 2 in \cite{GF}. By substituting (\ref{jkl}) in this formula we find
\[
\sum_{n=1}^\infty \phi_a(n)\frac{x^n}{1-x^n} = \sum_{n=1}^\infty x^n \tau(\gcd an) n.
\]
Multiplying the right hand side by $(1-2x^a+x^{2a})$ yields
\begin{align*}
&( \sum_{n=1}^\infty x^n \tau(\gcd an) n) -2 ( \sum_{n=a+1}^\infty x^n \tau(\gcd an) (n-a)) \\
&+( \sum_{n=1+2a}^\infty x^n \tau(\gcd an) (n-2a)) = \sum_{n=1}^\infty c[a](n) x^n,
\end{align*}
where
\[
c[a](n)=\left\{ \begin{array}{ll} 
 \tau(\gcd an)n & 0<n\leq a, \\
 \tau(\gcd an)(n-2(n-a))=\tau(\gcd an)(2a-n) & a<n\leq 2a, \\
 \tau(\gcd an)(n-2(n-a)+n-2a)=0 & n>2a. \end{array} \right.
 \]
Rewriting, using (\ref{ida}), the fact that $\gcd a{a+k}=\gcd a{a-k}$, and dividing by $x$,
yields the result.
 \hfill $\square$
 
\noindent
\begin{center}
\includegraphics[width=6cm]{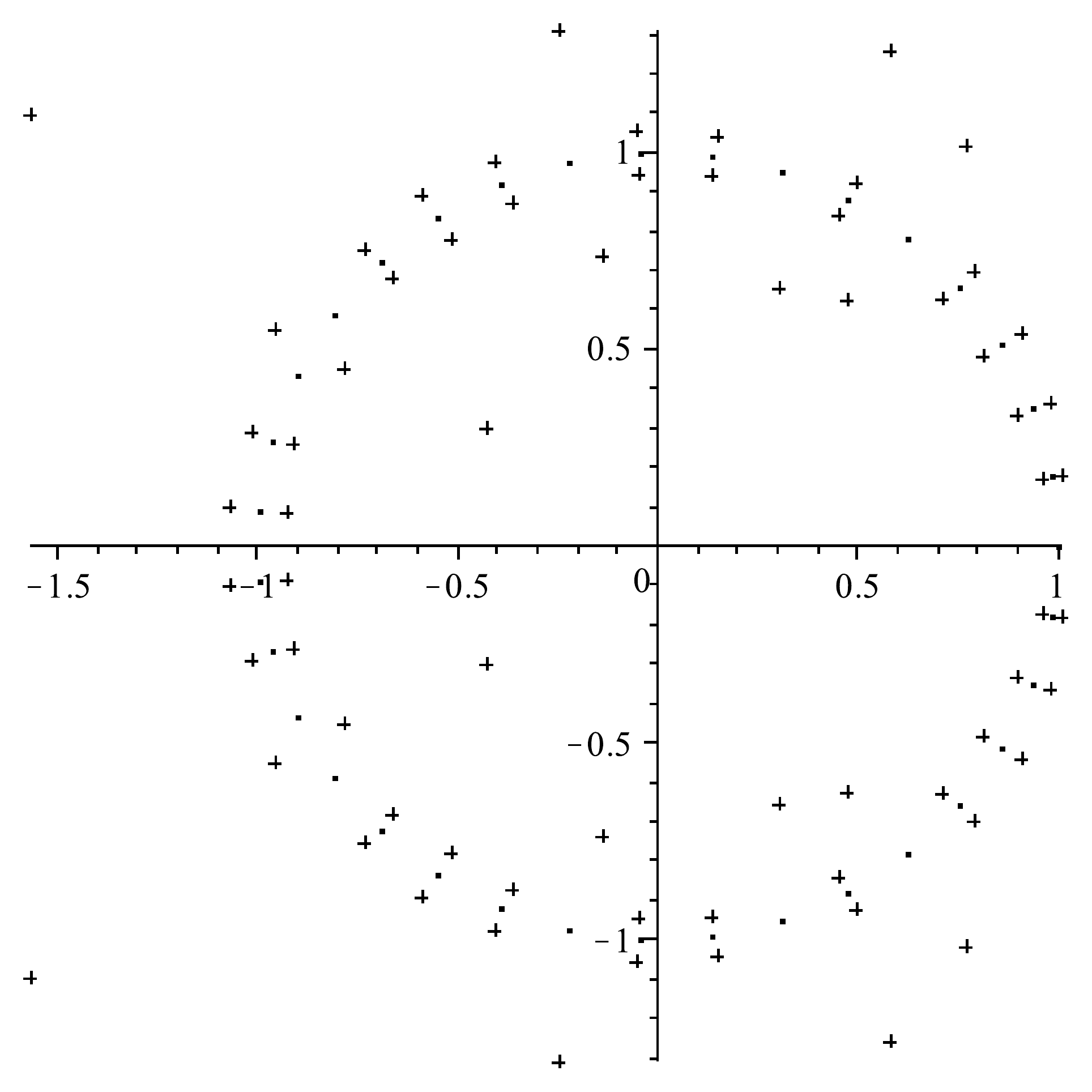}
\captionof{figure}{\label{35} The roots of $p[35]$ are depicted as crosses and the $35^{\text{th}}$ roots of unity as points.
This figure shows that roots of $p[a]$ are closest to primitive $a^{\text{th}}$ roots of unity.}
\end{center}

\subsection{A series related to the Lambert series of $\phi_a$}
Liouville \cite{Liou} also showed
\[
\sum_{m=1}^\infty \phi(m) \frac{x^m}{1+x^m} = (1+x^2)\frac{x}{(1-x^2)^2}.
\]
We show that
\begin{equation} \label{mola}
\sum_{m=1}^\infty \phi_a(m) \frac{x^m}{1+x^m} = q[a](x)\frac{x}{(1-x^{2a})^2},
\end{equation}
where $q[a](x)=\sum_{k=1}^{4a} b[a](k) x^{k-1}$, with
\[
b[a]=h[a]\circ t[2a],\quad h[a](k)=\pi_a(k)-2[2\mid k]\pi_a(k/2).
\]
At the end of this section we show that $1+x^2$ divides $q[a](x)$ if $a$ is odd.

\smallskip
\noindent
{\bf Proof} [of (\ref{mola})] As the left hand side of (\ref{mola}) is obtained from
the left hand side of (\ref{lam}) by subtracting twice the same series with $x$ replaced
by $x^2$, the same is true for the right hand side. Thus it follows that
$
q[a](x)=p[a](x)(1+x^a)^2-2p[a](x^2)x,
$
and hence, that
\[
b[a](k)=\left\{ \begin{array}{ll}
\pi_a(k)-2[2\mid k]\pi_a(k/2)& k\leq a, \\
\pi_a(2a-k)+2\pi_a(k-a)-2[2\mid k]\pi_a(k/2) & a<k\leq 2a, \\
2\pi_a(3a-k)+\pi_a(k-2a)-2[2\mid k]\pi_a(2a-k/2) & 2a<k\leq 3a, \\
\pi_a(4a-k)-2[2\mid k]\pi_a(2a-k/2) & 3a<k\leq 4a.
\end{array} \right.
\]
The result follows due to the identities
\[
\pi_a(2a-k)+2\pi_a(k-a)=\pi_a(k),\ \ 
2\pi_a(3a-k)+\pi_a(k-2a)=\pi_a(4a-k),
\]
which are easily verified using (\ref{ida}).  \hfill $\square$

\smallskip
\noindent
We can express the functions $b[a]$ and $h[a]$ in terms of an interesting fractal function.
Let a function $\kappa$ of two variables be defined recursively by
\begin{equation} \label{deka}
\kappa[a](n)=\left\{ \begin{array}{ll} 0 & 2\mid n , 2\nmid a, \text{ or } n=0, \\
\kappa[a/2](n/2) & 2\mid n , 2\mid a, \\
\tau(\gcd an) & 2\nmid n. \end{array} \right.
\end{equation}
The following properties are easily verified using the definition.
We have
\begin{equation} \label{pr1}
\kappa[a](2a+n)=\kappa[a](2a-n),
\end{equation}
and, with $\gcd ab=1$,
\begin{equation} \label{pr2}
\kappa[an](bn)=\left\{ \begin{array}{ll} 0 & 2 \mid b, \\ \alpha(n) & 2\nmid b, \end{array} \right.
\end{equation}
where $\alpha$ denotes the number of odd divisors function, i.e. for all $k$ and odd $m$
\begin{equation} \label{alp}
\alpha(2^k m) = \tau(m).
\end{equation}
Property (\ref{pr2}) is a quite remarkable fractal property; from the origin in every direction we see either the
zero sequence, or $\alpha$, at different scales.

We claim that
\begin{equation} \label{e1}
h[a]=\kappa[a]\pi
\end{equation}
follows from (\ref{deka}), (\ref{alp}), and (\ref{ida}). From (\ref{e1}) and (\ref{pr1}) we obtain
\begin{equation} \label{e2}
b[a]=\kappa[a]t[2a].
\end{equation}
We now prove that $1+x^2$ divides $q[a]$ when $a$ is odd.
Noting that $b[a](2a+k)=b[a](2a-k)$ and,
when $2\nmid a$, $b[a](2k)=0$, we therefore have
\begin{align*}
q[a](x)&=\sum_{n=1}^{2a} b[a](2n-1)x^{2n-2} \\
&=\sum_{m=1}^{a} b[a](2a-2m+1)x^{2a-2m} + b[a](2a+2m-1)x^{2a+2m-2} \\
&=\sum_{m=1}^{a} b[a](2a-2m+1)x^{2a-2m}(1+x^{4m-2}),
\end{align*}
which vanishes at the points where $x^2=-1$.  \hfill $\square$

\smallskip
\noindent
Apart from the factor $1+x^2$ when $a$ is odd, the polynomial $q[a]$ seems to be irreducible over ${\mathbb Z}$
and its zeros are in some sense close to the $2a^{\text{th}}$ roots of $-1$ or, to the $(a+1)^{\text{st}}$ roots of unity, see Figure \ref{19}.

\noindent
\begin{center}
\includegraphics[width=6cm]{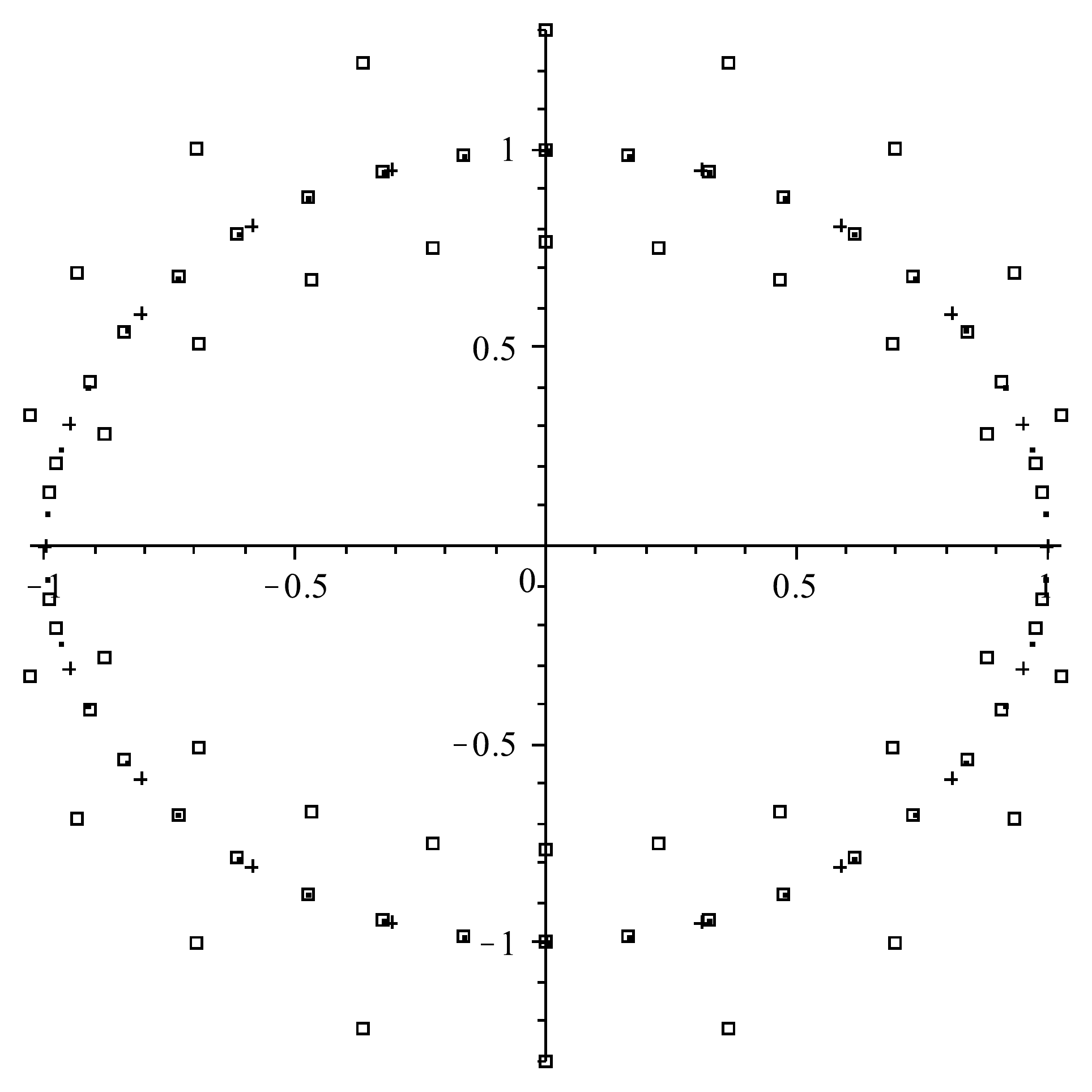}
\captionof{figure}{\label{19} The roots of $q[19]$ are depicted as boxes, the $38^{\text{th}}$ roots of $-1$ as points, and the
$20^{\text{th}}$ roots of unity as crosses. %This figure shows that, when $a$ is prime, the zeros of $q_a$ that are close to $2a$-th roots of -1 get closer when they are near $\pm i$, whereas the others get closer when they are close to $\pm 1$.
}
\end{center}

\subsection{A perfect square}\smallskip
Our last identity generalises the faint fact that $\phi(1)=1$. We have
\begin{equation} \label{nsq}
\sum_{a=1}^n \phi_a(n) = n^2.
\end{equation}

\smallskip
\noindent
{\bf Proof} [of (\ref{nsq})] For any lattice point $(i,j)$ in the square $[1,n]\times[1,n]$ the product $i\cdot j \mod n$ is congruent
to some $a$ in the range $[1,n]$. \hfill $\square$

\bigskip
\noindent
{\bf Acknowledgment} This research has been funded by the Australian Research Council through the Centre
of Excellence for Mathematics and Statistics of Complex Systems.


\begin{thebibliography}{10}
{%\small
\bibitem{AA} D.R. Anderson and T.M. Apostol, The evaluation of Ramanujan's sum and generalizations, Duke Math. J. 20 (1953) 211-216.

\bibitem{Tom} T.M. Apostol, Introduction to Analytic Number Theory, (1976) Springer-Verlag, New York.

\bibitem{Tom2} T.M. Apostol, Arithmetical properties of generalized Ramanujan sums, Pacific J. Math. 41 (1972) 281-293.

\bibitem{Bro} K.A. Broughan (2001), "The gcd-sum function", Journal of Integer Sequences 4, Art. 01.2.2.

\bibitem{Dick} L.E. Dickson, History of the theory of numbers, Carnegie Inst., Washington, D.C., 1919; reprinted by Chelsea, New York, 1952.

\bibitem{Dir} G.L. Dirichlet, \"Uber die Bestimmung der mittleren Werte in der Zahlentheorie, Abh. Akad. Wiss. Berlin, 1849; 78-8. Also in Werke, vol. 2, 1897, 60-64.

\bibitem{Hol} O. H\"older, Zur Theorie der Kreisteilungsgleichung $K_m(x)=0$, Prace Matematyczno
Fizyczne 43 (1936) 13-23.

\bibitem{kluyver} J.C. Kluyver, Some formulae concerning the integers less than $n$ and prime to $n$, in: KNAW, Proceedings 9 I, Amsterdam, (1906) 408-414.

\bibitem{Liou}
J. Liouville, Sur quelques s\'eries et produits infinis,
J. Math. Pures Appl. 2 (1857) 433-440.

\bibitem{LTPH}
L. T\'oth and P. Haukkanen, The discrete Fourier transform of r-even functions, Acta Univ. Sapientiae Math. 3 (2011) 5-25.

\bibitem{Ram} S. Ramanujan, On Certain Trigonometric Sums and their Applications in the Theory of Numbers, Transactions of the Cambridge Philosophical Society 22 (1918) 259Ð276. Also in: Collected papers of Srinivasa Ramanujan, Ed. G.H. Hardy et al, Chelsea Publ. Comp., New York (1962) 179-199.

\bibitem{SC}
J. Sandor and B. Crstici, Handbook of Number Theory II, (2004) Kluwer Acad. Publ., Dordrecht.

\bibitem{VoS}
R. Daublebsky Von Sterneck, Ein Analogon zur additiven Zahlentheorie, Sitzber, Akad. Wiss. Wien, Math. Naturw. Klasse, vol. I ll (Abt. IIa) (1902), 1567-1601.


\bibitem{GF} H. S. Wilf, Generatingfunctionology, Academic Press, 2nd edition, 1994.

\bibitem{WS} W. Schramm, The Fourier transform of functions of the greatest common divisor, Integers 8 (2008) A50.
}
\end{thebibliography}
\end{document}